\documentclass[a4paper,12pt]{article}
\usepackage{comment}
\usepackage{cite}
\usepackage{amsmath}
\usepackage{amssymb}
\usepackage{amsfonts}
\usepackage[T1]{fontenc}
\usepackage[utf8]{inputenc}
\usepackage{graphicx}
\usepackage{fancyhdr}
\usepackage{float}
\usepackage{xcolor}
\usepackage{authblk}
\usepackage{mathrsfs}
\usepackage{empheq}
\usepackage[hyphens]{url}
\usepackage{hyperref} 
\usepackage[]{breakurl}
\usepackage{amsmath}
\usepackage{amssymb}

\usepackage{geometry}
\begin{document}

\title{Integral representations of Catalan numbers using Touchard-like identities}

\author[$\dagger$]{Jean-Christophe {\sc Pain}$^{1,2,}$\footnote{jean-christophe.pain@cea.fr}\\
\small
$^1$CEA, DAM, DIF, F-91297 Arpajon, France\\
$^2$Universit\'e Paris-Saclay, CEA, Laboratoire Mati\`ere en Conditions Extr\^emes,\\ 
F-91680 Bruy\`eres-le-Ch\^atel, France
}

\date{}

\maketitle

\begin{abstract}
In this article, we use the Touchard identity in order to obtain new integral representations for Catalan numbers. The main idea consists in combining the identity with a known integral representation and resorting to the binomial theorem. The same procedure is applied to a variant of the Touchard identity proposed by Callan a few years ago. The method presented here can be generalized to derive additional integral representations from known ones, provided that the latter have a well-suited form and lend themselves to an analytical summation under the integral sign.
\end{abstract}

\section{Introduction}

The Catalan numbers \cite{Stanley1999,Koshy2008,Stanley2015}:
\begin{equation}
    C_n=\frac{1}{n+1}\binom{2n}{n}
\end{equation}
are still a subject of investigation in number theory. $C_n$ is the number of different ways $n + 1$ factors can be completely parenthesized, but also the number of Dyck words of length $2n$ or the number of full binary trees with $n + 1$ leaves (or $n$ internal nodes). Several integral representations already exist, ontained from Mellin transforms \cite{Penson2001}, or using recurrence relations \cite{Dana2005} for instance. Many interesting integral representations were reviewed in Ref. \cite{Qi2017}.

In the present work, we start from the known representation \cite{Guo2023}:
\begin{equation}\label{int0}
    C_n=\frac{2}{\pi}\int_0^{\infty}\frac{t^2}{\left(t^2+\frac{1}{4}\right)^{n+2}}\,\mathrm{d}t.
\end{equation}
Starting from the latter integral form, we obtain new ones, which can not be easily deduced from the integral form (\ref{int0}).

\section{The Touchard identity for Catalan numbers}

In 1928, Touchard published the following identity \cite{Touchard1928}:
\begin{equation}\label{tou}
C_{n+1}=\sum_{k=0}^{\lfloor\frac{n}{2}\rfloor}\binom{n}{2k}\,C_k\,2^{n-2k}.    
\end{equation}
Different proofs of identity (\ref{tou}) were proposed (see for instance Ref. \cite{Shapiro1976}), even rather recently \cite{Regev2015}. Let us introduce the expression (\ref{int0}) in Eq. (\ref{tou}). This yields
\begin{equation}
    C_{n+1}=\frac{2}{\pi}\int_0^{\infty}\frac{t^2}{\left(t^2+\frac{1}{4}\right)^2}\sum_{k=0}^{\lfloor\frac{n}{2}\rfloor}\binom{n}{2k}\frac{2^{n-2k}}{\left(t^2+\frac{1}{4}\right)^k}\,\mathrm{d}t.
\end{equation}
Then one has, using the binomial theorem, the equality (for $a>0$):
\begin{equation}
\sum_{k=0}^{\lfloor\frac{n}{2}\rfloor}\binom{n}{2k}\,a^k\,b^{n-2k}=\frac{1}{2}\left[(b-\sqrt{a})^n+(b+\sqrt{a})^n\right],
\end{equation}
yielding, with $a=1/(t^2+1/4)^k$ and $b=2$:
\begin{equation}\label{new1}
    C_{n}=\frac{1}{\pi}\int_0^{\infty}\frac{t^2}{\left(t^2+\frac{1}{4}\right)^2}\left[\left(2-\frac{1}{\sqrt{t^2+\frac{1}{4}}}\right)^{n-1}+\left(2+\frac{1}{\sqrt{t^2+\frac{1}{4}}}\right)^{n-1}\right]\,\mathrm{d}t,
\end{equation}
which is the first main result of the present work. Equation (\ref{int0}) can be put in the form 
\begin{equation}\label{int0bis}
    C_n=\frac{2}{\pi}\int_0^{\infty}\frac{t^2}{\left(t^2+\frac{1}{4}\right)^2}\,f_1(t)\,\mathrm{d}t,
\end{equation}
with
\begin{equation}
    f_1(t)=\frac{1}{\left(t^2+\frac{1}{4}\right)^n}
\end{equation}
and the new integral can be written as (\ref{new1}):
\begin{equation}\label{new1bis}
    C_n=\frac{2}{\pi}\int_0^{\infty}\frac{t^2}{\left(t^2+\frac{1}{4}\right)^2}\,f_2(t)\,\mathrm{d}t,
\end{equation}
with
\begin{equation}
    f_2(t)=\frac{1}{2}\left[\left(2-\frac{1}{\sqrt{t^2+\frac{1}{4}}}\right)^{n-1}+\left(2+\frac{1}{\sqrt{t^2+\frac{1}{4}}}\right)^{n-1}\right]\,\mathrm{d}t.
\end{equation}
It is worth noting that integral (\ref{new1}) is different from Eq. (\ref{int0}), since, as can be seen in Fig. \ref{fig1}, $f_1$ and $f_2$ are different (they intersect only for one value).

\section{Callan's variant of Touchard's identity}

The Catalan number $C_n$ counts the number of dissections of a regular polygon with $n+2$ sides into triangles. Counting the dissections by numbering the triangles containing two sides of the polygon among their three edges, Callan obtained, in 2013, the following relation, valid for $n>1$ \cite{Callan2013}:
\begin{equation}\label{cal}
    C_n=\sum_{k=1}^{\lfloor\frac{n}{2}\rfloor}2^{n-2k}\binom{n}{2k}\,C_k\,\frac{k(n+2)}{n(n-1)}.
\end{equation}
Callan also illustrated the connection of Eq. (\ref{cal}) with Touchard's identity.

Inserting expression (\ref{int0}) into Eq. (\ref{cal}) yields, for $n>1$:
\begin{equation}
    C_{n}=\frac{2}{\pi}\frac{(n+2)}{n(n+1)}\int_0^{\infty}\frac{t^2}{\left(t^2+\frac{1}{4}\right)^2}\sum_{k=0}^{\lfloor\frac{n}{2}\rfloor}k\binom{n}{2k}\frac{2^{n-2k}}{\left(t^2+\frac{1}{4}\right)^k}\,\mathrm{d}t.
\end{equation}
One can show, using the binomial theorem, that (still with $a>0$):
\begin{equation}
    \sum_{k=0}^{\lfloor\frac{n}{2}\rfloor}\binom{n}{2k}\,k\,a^k\,b^{n-2k}=\frac{n\sqrt{a}}{4}\left[(b+\sqrt{a})^{n-1}-(b-\sqrt{a})^{n-1}\right],
\end{equation}
leading to
\begin{equation}\label{new2}
    C_{n}=\frac{(n+2)}{2(n-1)\pi}\int_0^{\infty}\frac{t^2}{\left(t^2+\frac{1}{4}\right)^{5/2}}\left[\left(2+\frac{1}{\sqrt{t^2+\frac{1}{4}}}\right)^{n-1}-\left(2-\frac{1}{\sqrt{t^2+\frac{1}{4}}}\right)^{n-1}\right]\,\mathrm{d}t,
\end{equation}
which constitutes the second main result of the present work. The latter integral (\ref{new2}) can be put in the form
\begin{equation}\label{new2bis}
    C_n=\frac{2}{\pi}\int_0^{\infty}\frac{t^2}{\left(t^2+\frac{1}{4}\right)^2}\,f_3(t)\,\mathrm{d}t,
\end{equation}
with
\begin{equation}
    f_3(t)=\frac{(n+2)}{4(n-1)}\frac{1}{\sqrt{t^2+\frac{1}{4}}}\left[\left(2+\frac{1}{\sqrt{t^2+\frac{1}{4}}}\right)^{n-1}-\left(2-\frac{1}{\sqrt{t^2+\frac{1}{4}}}\right)^{n-1}\right]\,\mathrm{d}t.
\end{equation}

\begin{figure}[!ht]
    \centering
    \includegraphics[width=11cm]{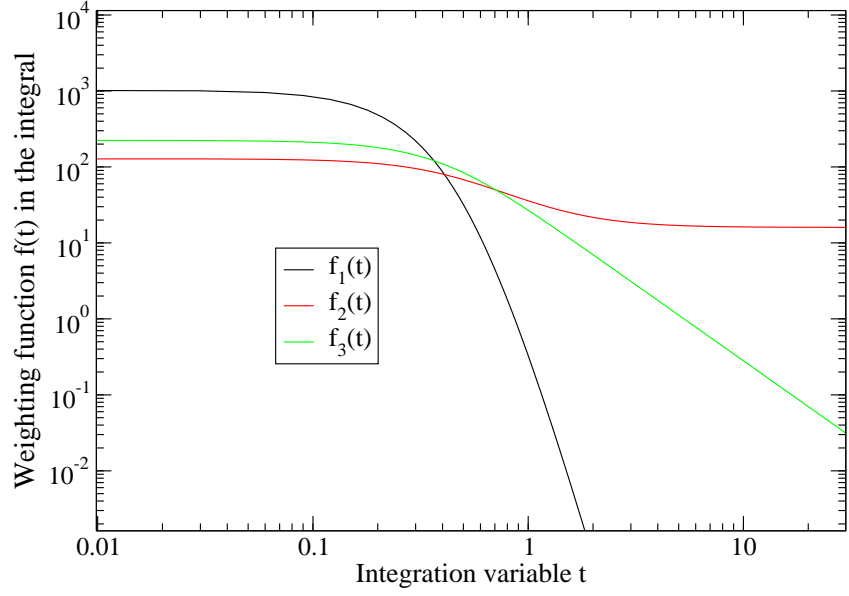}
    \caption{Functions $f_1$, $f_2$ and $f_3$ involved as weighting functions in the three integral representations ((\ref{int0})-(\ref{int0bis})), ((\ref{new1})-(\ref{new1bis})) and ((\ref{new2})-(\ref{new2bis})).}\label{fig1}
\vspace{1cm}
\end{figure}
Here also the latter integral is different from Eq. (\ref{int0}) and (\ref{new1}), since, as can be seen in Fig. \ref{fig1}, $f_3$ is different from $f_1$ and $f_2$.

\begin{figure}[!ht]
\vspace{1cm}
    \centering
    \includegraphics[width=11cm]{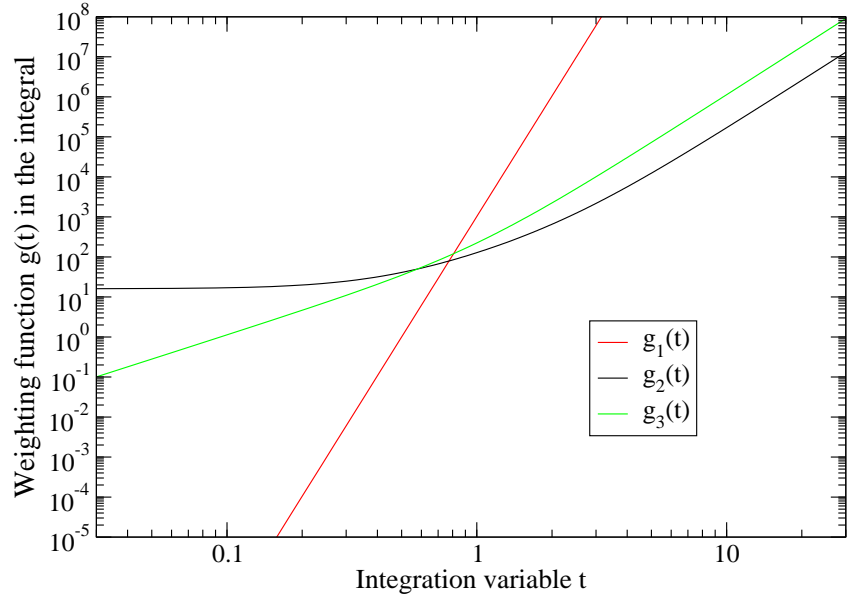}
    \caption{Functions $g_1$, $g_2$ and $g_3$ involved as weighting functions in the integral representations (\ref{bas}), (\ref{bis}) and (\ref{ter}).}\label{fig2}
\end{figure}

Other integral representations can be derived in the same way. For instance, using the well-known integral representation \cite{Choi2020}:
\begin{equation}\label{bas}
    C_n=\frac{2^{2n+1}}{\pi}\int_{-1}^1t^{2n}\,\sqrt{1-t^2}\,\mathrm{d}t,
\end{equation}
together with:
\begin{equation}
    \sum_{k=0}^{\lfloor\frac{n}{2}\rfloor}\binom{n}{2k}\,a^{2k}\,b^{n-2k}=\frac{1}{2}\left[(a+b)^{n}+(b-a)^{n}\right],    
\end{equation}
one gets, using the Touchard identity (\ref{tou}), the representation
\begin{equation}\label{bis}
C_n=\frac{2^{n-1}}{\pi}\int_{-1}^1\sqrt{1-t^{2}}\,\left[(1-t)^{n-1}+(1+t)^{n-1}\right]\,\mathrm{d}t,
\end{equation}
and using Callan's variant (\ref{cal}), for $n>1$, combined with
\begin{equation}
    \sum_{k=0}^{\lfloor\frac{n}{2}\rfloor}\binom{n}{2k}\,k\,a^{2k}\,b^{n-2k}=\frac{n\,a}{4}\left[(a+b)^{n-1}-(b-a)^{n-1}\right],    
\end{equation}
one obtains the integral representation
\begin{equation}\label{ter}
C_n=\frac{2^{n-1}(n+2)}{\pi(n-1)}\int_{-1}^1\,t\,\sqrt{1-t^{2}}\,\left[(t+1)^{n-1}-(1-t)^{n-1}\right]\,\mathrm{d}t.
\end{equation}
In that case also the representations (\ref{bis}) and (\ref{ter}) are new in the sense that they cannot be easily deduced from Eq. (\ref{bas}), since the functions 
\begin{equation}
g_1(t)=(2t)^{2n},
\end{equation}
\begin{equation}
g_2(t)=2^{n-2}\,\left[(1-t)^{n-1}+(1+t)^{n-1}\right]
\end{equation}
and for $n>1$:
\begin{equation}
g_3(t)=2^{n-2}\,\frac{(n+2)}{(n-1)}\,t\,\left[(t+1)^{n-1}-(1-t)^{n-1}\right]
\end{equation}
are clearly different, as displayed in Fig. \ref{fig2}.


\section{Conclusion}

We derived new integral representations for Catalan numbers, inserting a well-chosen known integral representation of the numbers in the Touchard identity. The same procedure was applied to a variant of the Touchard identity published by Callan. The method presented here can be generalized to derive additional integral representations from known ones, provided that the latter have an appropriate form, yielding a simplification of the finite sum. It is also possible to extend the technique described in this paper to other important numbers, such as the $q$-Narayana numbers for instance \cite{Pan2022}.

\end{document}